# On The Exact Solution of Newell-Whitehead-Segel Equation Using the Homotopy Perturbation Method

[1]S. S. Nourazar, [1] M. Soori and [2]A. Nazari-Golshan

[1]Mechanical Engineering Department, Amirkabir University of Technology, Tehran, Iran
[2]Department of Physics, and cAmirkabir University of Technology, Tehran, Iran

**Abstract:** In the present work, we use the homotopy perturbation method (HPM) to solve the Newell-Whitehead-Segel non-linear differential equations. Four case study problems of Newell-Whitehead-Segel are solved by the HPM and the exact solutions are obtained. The trend of the rapid convergence of the sequences constructed by the method toward the exact solution is shown numerically. As a result the rapid convergence towards the exact solutions of HPM indicates that, using the HPM to solve the Newell-Whitehead-Segel non-linear differential equations, a reasonable less amount of computational work with acceptable accuracy may be sufficient. Moreover the application of the HPM proves that the method is an effective and simple tool for solving the Newell-Whitehead-Segel non-linear differential equations.

## INTRODUCTION

Recently lots of attentions are devoted toward the semi-analytical solution of real-life mathematical modeling that is inherently nonlinear differential equations. Most of the nonlinear differential equations do not have an analytical solution. The idea of the homotopy perturbation method was first pioneered by (He, 1999) Later the homotopy perturbation method (HPM) which is a semi-analytical method is applied to solve the non-linear non-homogeneous partial differential equations (He, 2005; Yildirim, 2009; He, 2006; He, 2006; He, 2006; He, 2005; He, 1998; He, 1999; He, 1999; He, 2000; He, 2006; He, 2005; He, 2005; Koçak et al., 2011; Gepreel, 2011; Cao and Bo. Han, 2011; Wazwaz, 2009). Ezzati and Shakibi (Ezzati, 2011) solved the Newell-Whitehead equation using the Adomian decomposition and multi-quadric quasi-interpolation methods. They concluded that the Adomian decomposition and multi-quadric quasi-interpolation methods are reasonable methods to solve the Newell-Whitehead equation with acceptable accuracy.

In the present work, the homotopy perturbation method (HPM) is applied to obtain the closed form solution of the non-linear Newell-Whitehead-Segel equation. Four case study problems of non-linear Newell-Whitehead-Segel equations are solved using the HPM. The trend of the rapid convergence towards the exact solution is shown when compared to the exact solution. The Newell-Whitehead-Segel equation models the interaction of the effect of the diffusion term with the nonlinear effect of the reaction term. The Newell-Whitehead-Segel equation is written as:

$$\frac{\partial u}{\partial t} = k \frac{\partial^2 u}{\partial x^2} + au - bu^q \tag{1}$$

Where $a$, $b$ and $k$ are real numbers with $k > 0$, and $q$ is a positive integer. In Eq. (1) the first term on the left hand side, $\frac{\partial u}{\partial t}$, expresses the variations of $u(x,t)$ with time at a fixed location, the first term on the right hand side, $\frac{\partial^2 u}{\partial x^2}$, expresses the variations of $u(x,t)$ with spatial variable $x$ at a specific time and the remaining terms on the right hand side, $au - bu^q$, takes into account the effect of the source term. In Eq. (1) $u(x,t)$ is a function of the spatial variable $x$ and the temporal variable $t$, with $x \in R$ and $t \geq 0$. The function $u(x,t)$ may be thought of as the (nonlinear) distribution of temperature in an infinitely thin and long rod or as the flow velocity of a fluid in an infinitely long pipe with small diameter. The Newell-Whitehead-Segel equations have wide applicability in mechanical and chemical engineering, ecology, biology and bio-engineering.

**Corresponding Author:** S.S. Nourazar, Mechanical Engineering Department, Amirkabir University of Technology, Tehran, Iran
E-mail: icp@aut.ac.ir

### 2. The Idea Of Homotopy Perturbation Method:

The homotopy perturbation method (HPM) is originally initiated by (He, 1999; He, 2005; Yildirim, 2009; He, 2006; He, 2006; He, 2006; He, 2005; He, 1998; He, 1998; He, 1999; He, 1999; He, 2000; He, 2006; He, 2005; He, 2005; He, 2005; Koçak et al., 2011). This is a combination of the classical perturbation technique and homotopy techinique. The basic idea of the HPM for solving nonlinear differential equations is as follow; consider the following differential equation:

$$E(u) = 0, \tag{2}$$

$$H(u, p) = (1 - p)(F(u) - H(u)) + p\big(E(u)\big) = 0. \tag{3}$$

Where $E(u)$ is any differential operator. We construct a homotopy as follow:

Where are *F(u)*, *H(u)* are functional operators with the known solution v₀. It is clear that when *p* is equal to zero then $H(u, 0) = F(u) - H(u) = 0$, and when p is equal to 1, $H(u, 1) = E(u) = 0$. It is worth nothing that is the embedding parameter *p* increases monotonically from zero to unity the zero order solution v₀ continuously deforms into the original problem $E(u) = 0$. The embedding parameter, $p \in [0,1]$, is considered as an expending parameter(He, 1999). In the homotopy perturbation method the embedding parameter *p* is used to get series expansion for solution as:

$$u = \sum_{i=0}^{\infty} p^i v_i = v_0 + p v_1 + p^2 v_2 + p^3 v_3 + \cdots \tag{4}$$

When $p \longrightarrow 1$, then Eq. (3) becomes the approximate solution to Eq. (2) as:

$$u = v_0 + v_1 + v_2 + v_3 + \cdots \tag{5}$$

The series Eq. (5) is a convergent series and the rate of convergence depends on the nature of Eq. (2) (He, 1999; He, 2005; Yildirim, 2009; He, 2006; He, 2006; He, 2006; He, 2005; He, 1998; He, 1998; He, 1999; He, 1999; He, 2000; He, 2006; He, 2005; He, 2005; Koçak et al., 2011). It is also assumed that Eq. (3) has a unique solution and by comparing the like powers of *p* the solution of various orders is obtained. These solutions are obtained using the Maple package.

### 3. The Newell-Whitehead-Segel Equation:

To illustrate the capability and reliability of the method, four cases of nonlinear diffusion equations are presented.

***Case I:***

In Eq. (1) for $a = 2, b = 3, k = 1$ and $q = 2$ the Newell-Whitehead-Segel equation is written as:

$$\frac{\partial u}{\partial t} = \frac{\partial^2 u}{\partial x^2} + 2u - 3u^2 \tag{6}$$

Subject to a constant initial condition

$$u(x, 0) = \lambda, \tag{7}$$

We construct a homotopy for Eq. (6) in the following form:

$$H(v, p) = (1 - p)\left[\frac{\partial v}{\partial t} - \frac{\partial u_0}{\partial t}\right] + p\left[\frac{\partial v}{\partial t} - \frac{\partial^2 v}{\partial x^2} - 2v + 3v^2\right] = 0 \tag{8}$$

The solution of Eq. (6) can be written as a power series in $\boldsymbol{p}$ as:

$$v = v_0 + p\, v_1 + p^2 v_2 + \cdots \tag{9}$$

Substituting Eq. (9) and Eq. (7) into Eq. (8) and equating the terms with identical powers of $p$:

$$p^0 : \frac{\partial v_0}{\partial t} = \frac{\partial u_0}{\partial t}, \qquad\qquad v_0(x,0) = \lambda,$$

$$p^1 : \frac{\partial v_1}{\partial t} - \frac{\partial u_0}{\partial t} = \frac{\partial^2 v_0}{\partial x^2} + 2v_0 - 3v_0^2, \qquad v_1(x,0) = 0,$$

$$p^2 : \frac{\partial v_2}{\partial t} = \frac{\partial^2 v_1}{\partial x^2} - 3v_0 v_1 + v_1(2 - 3v_0), \qquad v_2(x,0) = 0, \qquad (10)$$

$$p^3 : \frac{\partial v_3}{\partial t} = \frac{\partial^2 v_2}{\partial x^2} - 3v_0 v_2 - 3v_1^2 + v_2(2 - 3v_0), \qquad v_3(x,0) = 0.$$

Using the Maple package to solve recursive sequences, Eq. (10), we obtain the followings:

$$v_0(x,t) = \lambda,$$

$$v_1(x,t) = \lambda(2 - 3\lambda)t,$$

$$v_2(x,t) = 2\lambda(2 - 3\lambda)(1 - 3\lambda)\frac{t^2}{2!}, \qquad\qquad (11)$$

$$v_3(x,t) = 2\lambda(2 - 3\lambda)(27\lambda^2 - 18\lambda + 2)\frac{t^3}{3!}.$$

By setting $p = 1$ in Eq. (9) the solution of Eq. (6) can be obtained as
$v = v_0 + v_1 + v_2 + v_3 + \dots$ Therefore the solution of Eq. (6) is written as:

$$v(x,t) = \lambda + \lambda(2 - 3\lambda)t + 2\lambda(2 - 3\lambda)(1 - 3\lambda)\frac{t^2}{2!} + 2\lambda(2 - 3\lambda)(27\lambda^2 - 18\lambda + 2)\frac{t^3}{3!} + \dots$$

$$(12)$$

The Taylor series expansion for $\dfrac{-\frac{2}{3}\lambda e^{2t}}{-\frac{2}{3} + \lambda - \lambda e^{2t}}$ is written as:

$$\frac{-\frac{2}{3}\lambda e^{2t}}{-\frac{2}{3} + \lambda - \lambda e^{2t}} = \lambda + \lambda(2 - 3\lambda)t + 2\lambda(2 - 3\lambda)(1 - 3\lambda)\frac{t^2}{2!}$$

$$+ 2\lambda(2 - 3\lambda)(27\lambda^2 - 18\lambda + 2)\frac{t^3}{3!} + \cdots \qquad (13)$$

By substituting Eq. (13) into Eq. (12), thus Eq. (12) can be rewritten as:

$$v(x,t) = \frac{-\frac{2}{3}\lambda e^{2t}}{-\frac{2}{3} + \lambda - \lambda e^{2t}} \qquad\qquad (14)$$

This is the exact solution of the problem, Eq.(6). Table 1 shows the trend of rapid convergence of the results of $S_0(x,t) = v_0(x,t)$ to $S_5(x,t) = \sum_{i=0}^{5} v_i(x,t)$ using the HPM. The rapid convergence of the solution toward the exact solution, the maximum relative error of less than 0.0046%, is achieved as shown in table 1.

Table 1 shows: he percentage of relative errors of the results of $S_0(x,t) = v_0(x,t)$ to $S_5(x,t) = \sum_{i=0}^{5} v_i(x,t)$ of the HPM solution of Eq. (6) for $\lambda$=0.1.

**Case II**: In Eq. (1) for $a = 1, b = 1, k = 1$ and $q = 2$ the Newell-Whitehead-Segel equation is written as:

**Table 1:**

| | | percentage of relative error (%RE) |
|---|---|---|
| $t = 0.1$ | $S_0(x,t)$ | 0.1540788596 |
| | $S_1(x,t)$ | 0.01027226573 |
| | $S_2(x,t)$ | 0.0002058041542 |
| | $S_3(x,t)$ | 0.00001949256084 |
| | $S_4(x,t)$ | 0.000001708760704 |
| | $S_5(x,t)$ | 3.214500334 e-8 |
| $t = .3$ | $S_0(x,t)$ | 0.3835101093 |
| | $S_1(x,t)$ | 0.06910026503 |
| | $S_2(x,t)$ | 0.003074197736 |
| | $S_3(x,t)$ | 0.001358981068 |
| | $S_4(x,t)$ | 0.0003091665977 |
| | $S_5(x,t)$ | 0.00001230390524 |
| $t = .4$ | $S_0(x,t)$ | 0.4680703804 |
| | $S_1(x,t)$ | 0.1063582391 |
| | $S_2(x,t)$ | 0.005078839539 |
| | $S_3(x,t)$ | 0.003988077959 |
| | $S_4(x,t)$ | 0.001125247109 |
| | $S_5(x,t)$ | 0.00004585924829 |

$$\frac{\partial u}{\partial t} = \frac{\partial^2 u}{\partial x^2} + u - u^2 \tag{15}$$

Subject to initial condition

$$u(x,0) = \frac{1}{\left(1 + e^{\frac{x}{\sqrt{6}}}\right)^2} \tag{16}$$

To solve Eq. (15) we construct a homotopy in the following form:

$$H(v,p) = (1-p)\left[\frac{\partial v}{\partial t} - \frac{\partial u_0}{\partial t}\right] + p\left[\frac{\partial v}{\partial t} - \frac{\partial^2 v}{\partial x^2} - v(1-v)\right] = 0. \tag{17}$$

The solution of Eq. (15) can be written as a power series in $p$ as:

$$v = v_0 + p\,v_1 + p^2 v_2 + \cdots \tag{18}$$

Substituting Eq. (18) and Eq. (16) in to Eq. (17) and equating the term with identical powers of $p$, leads to

$$p^0 : \frac{\partial v_0}{\partial t} = \frac{\partial u_0}{\partial t}, \qquad\qquad v_0(x,0) = \frac{1}{\left(1 + e^{\frac{x}{\sqrt{6}}}\right)^2}$$

$$p^1 : \frac{\partial v_1}{\partial t} + \frac{\partial u_0}{\partial t} = \frac{\partial^2 v_0}{\partial x^2} + v_0(1 - v_0), \qquad v_1(x,0) = 0,$$

$$p^2 : \frac{\partial v_2}{\partial t} = \frac{\partial^2 v_1}{\partial x^2} - v_0 v_1 + v_1(1 - v_0), \qquad v_2(x,0) = 0, \tag{19}$$

$$p^3 : \frac{\partial v_3}{\partial t} = \frac{\partial^2 v_2}{\partial x^2} + v_2 - v_1^2 - 2v_0 v_2, \qquad v_3(x,0) = 0.$$

Using the Maple package to solve recursive sequences, Eq. (19), we obtain the followings:

$$v_0(x,t) = \frac{1}{\left(1+e^{\frac{x}{\sqrt{6}}}\right)^2}$$

$$v_1(x,t) = \frac{5}{3}\frac{e^{\frac{x}{\sqrt{6}}}}{\left(1+e^{\frac{x}{\sqrt{6}}}\right)^3}t \qquad (20)$$

$$v_2(x,t) = \frac{25}{18}\left(\frac{e^{\frac{x}{\sqrt{6}}}\left(-1+2e^{\frac{x}{\sqrt{6}}}\right)}{\left(1+e^{\frac{x}{\sqrt{6}}}\right)^4}\right)\frac{t^2}{2}$$

$$v_3(x,t) = \left(\frac{125}{216}\frac{e^{\frac{x}{\sqrt{6}}}\left(4(e^{\frac{x}{\sqrt{6}}})^2-7e^{\frac{x}{\sqrt{6}}}+1\right)}{\left(1+e^{\frac{x}{\sqrt{6}}}\right)^5}\right)\frac{t^3}{3}$$

By setting $p = 1$ in Eq. (18) the solution of Eq. (15) can be obtained as $v = v_0 + v_1 + v_2 + v_3 + \dots$ Therefore the solution of Eq. (15) is written as:

$$v(x,t) = \frac{1}{\left(1+e^{\frac{x}{\sqrt{6}}}\right)^2} + \frac{5}{3}\frac{e^{\frac{x}{\sqrt{6}}}}{\left(1+e^{\frac{x}{\sqrt{6}}}\right)^3}t + \frac{25}{18}\left(\frac{e^{\frac{x}{\sqrt{6}}}\left(-1+2e^{\frac{x}{\sqrt{6}}}\right)}{\left(1+e^{\frac{x}{\sqrt{6}}}\right)^4}\right)\frac{t^2}{2} + \left(\frac{125}{216}\frac{e^{\frac{x}{\sqrt{6}}}\left(4(e^{\frac{x}{\sqrt{6}}})^2-7e^{\frac{x}{\sqrt{6}}}+1\right)}{\left(1+e^{\frac{x}{\sqrt{6}}}\right)^5}\right)\frac{t^3}{3} +$$

$$\dots$$

The Taylor series expansion for $\dfrac{1}{\left(1+e^{\frac{x}{\sqrt{6}}-\frac{5}{6}t}\right)^2}$ is written as $\qquad\qquad\qquad (21)$

$$\frac{1}{\left(1+e^{\frac{x}{\sqrt{6}}-\frac{5}{6}t}\right)^2} = \frac{1}{\left(1+e^{\frac{x}{\sqrt{6}}}\right)^2} + \frac{5}{3}\frac{e^{\frac{x}{\sqrt{6}}}}{\left(1+e^{\frac{x}{\sqrt{6}}}\right)^3}t + \frac{25}{18}\left(\frac{e^{\frac{x}{\sqrt{6}}}\left(-1+2e^{\frac{x}{\sqrt{6}}}\right)}{\left(1+e^{\frac{x}{\sqrt{6}}}\right)^4}\right)\frac{t^2}{2}$$

$$+ \left(\frac{125}{216}\frac{e^{\frac{x}{\sqrt{6}}}\left(4(e^{\frac{x}{\sqrt{6}}})^2-7e^{\frac{x}{\sqrt{6}}}+1\right)}{\left(1+e^{\frac{x}{\sqrt{6}}}\right)^5}\right)\frac{t^3}{3} + \dots \qquad (22)$$

By substituting Eq. (22) into Eq. (21), the Eq. (21) can be reduced to

$$v(x,t) = \frac{1}{\left(1+e^{\frac{x}{\sqrt{6}}-\frac{5}{6}t}\right)^2} \qquad (23)$$

This is the exact solution of the problem, Eq. (15). Table 2 shows the trend of rapid convergence of the results of $S_0(x,t) = v_0(x,t)$ to $S_5(x,t) = \sum_{i=0}^{5} v_i(x,t)$ using the HPM solution toward the exact solution. The maximum relative error of less than 0.0008% is achieved in comparison to the exact solution as shown in table 2.

**Table 2:**

| | | percentage of relative error (%RE) | | |
|---|---|---|---|---|
| | | $x = 1$ | $x = 1.5$ | $x = 1.8$ |
| $t = 0.1$ | $S_0(x,t)$ | 0.09374691950 | 0.1010110875 | 0.1051585828 |
| | $S_1(x,t)$ | 0.003020748923 | 0.003848096613 | 0.004359470571 |
| | $S_2(x,t)$ | 0.00001102987604 | 0.00002050070585 | 0.00004357185889 |
| | $S_3(x,t)$ | 0.000003008522247 | 0.000003116144165 | 0.000002949567767 |
| | $S_4(x,t)$ | 3.247853526 e-8 | 7.597179054 e-8 | 1.007266564 e-7 |
| | $S_5(x,t)$ | 2.475184340 e-9 | 6.284562692 e-10 | 4.103584300 e-10 |
| $t = .3$ | $S_0(x,t)$ | 0.2480807805 | 0.2663112457 | 0.2766522515 |
| | $S_1(x,t)$ | 0.02225395729 | 0.02841934500 | 0.03220845585 |
| | $S_2(x,t)$ | 0.0003852758255 | 0.0003051199126 | 0.0008095286291 |
| | $S_3(x,t)$ | 0.0002055820118 | 0.0002152872468 | 0.0002058262157 |
| | $S_4(x,t)$ | 0.000005574634536 | 0.00001431284116 | 0.00001929380478 |
| | $S_5(x,t)$ | 0.000001472642781 | 0.000001037290051 | 5.725167943 e-7 |
| $t = .4$ | $S_0(x,t)$ | 0.3115489693 | 0.3338572016 | 0.3464690907 |
| | $S_1(x,t)$ | 0.03586206901 | 0.04586954916 | 0.05200200785 |
| | $S_2(x,t)$ | 0.0009882331330 | 0.0004901161293 | 0.001569416719 |
| | $S_3(x,t)$ | 0.0005982451688 | 0.0006298759390 | 0.0006050514464 |
| | $S_4(x,t)$ | 0.00001947940975 | 0.00005317431566 | 0.00007241822783 |
| | $S_5(x,t)$ | 0.000007711072685 | 0.000005555853236 | 0.000003218048969 |

**Table 2** shows: the percentage of relative errors of the results of $S_0(x,t) = v_0(x,t)$ to $S_5(x,t) = \sum_{i=0}^{5} v_i(x,t)$ of the HPM solution of Eq. (15).

**Case III:** In Eq. (1) for $a = 1, b = 1, k = 1$ and $q = 4$ the Newell-Whitehead-Segel equation becomes:

$$\frac{\partial u}{\partial t} = \frac{\partial^2 u}{\partial x^2} + u - u^4 \tag{24}$$

Subject to initial condition

$$u(x,0) = \left(\frac{1}{1 + e^{\frac{5x}{\sqrt{10}}}}\right)^{\frac{2}{3}} \tag{25}$$

We construct a homotopy for Eq. (24) in the following form:

$$H(v,p) = (1-p)\left[\frac{\partial v}{\partial t} - \frac{\partial u_0}{\partial t}\right] + p\left[\frac{\partial v}{\partial t} - \frac{\partial^2 v}{\partial x^2} - v + v^4\right] = 0. \tag{26}$$

The solution of Eq. (24) can be written as a power series in $p$ as:

$$v = v_0 + p\,v_1 + p^2 v_2 + \cdots \tag{27}$$

Substituting Eq. (27) and Eq. (25) into Eq. (26) and equating the terms with identical powers of *p*:

$$p^0: \frac{\partial v_0}{\partial t} = \frac{\partial u_0}{\partial t}, \qquad v_0(x,0) = \frac{1}{\left(1 + e^{\frac{3}{\sqrt{10}}x}\right)^{\frac{2}{3}}}$$

$$p^1: \frac{\partial v_1}{\partial t} = \frac{\partial^2 u_0}{\partial x^2} - \frac{\partial u_0}{\partial t} + v_0(1 - v_0{}^3), \qquad v_1(x,0) = 0,$$

$$p^2: \frac{\partial v_2}{\partial t} = \frac{\partial^2 v_1}{\partial x^2} - 4v_0{}^3 v_1 + v_1, \qquad v_2(x,0) = 0, \qquad (28)$$

$$p^3: \frac{\partial v_3}{\partial t} = \frac{\partial^2 v_2}{\partial x^2} - 6v_0{}^2 v_1{}^2 - 4v_0{}^3 v_2 + v_2, \qquad v_3(x,0) = 0.$$

Using the Maple package to solve recursive sequences, Eq. (28), we obtain the followings:

$$v_0(x,t) = \left(\frac{1}{1 + e^{\frac{3}{\sqrt{10}}x}}\right)^{\frac{2}{3}}$$

$$v_1(x,t) = \frac{7}{5}\left(\frac{e^{\frac{3}{\sqrt{10}}x}}{\left(1 + e^{\frac{3}{\sqrt{10}}x}\right)^{\frac{5}{3}}}\right) t$$

$$v_2(x,t) = \frac{49}{50}\left(\frac{\left(2e^{\frac{3}{\sqrt{10}}x} - 3\right)e^{\frac{3}{\sqrt{10}}x}}{\left(1 + e^{\frac{3}{\sqrt{10}}x}\right)^{\frac{8}{3}}}\right)\frac{t^2}{2} \qquad (29)$$

$$v_3(x,t) = \frac{343}{1000}\left(\frac{\left(4\left(e^{\frac{3}{\sqrt{10}}x}\right)^2 - 27e^{\frac{3}{\sqrt{10}}x} + 9\right)e^{\frac{3}{\sqrt{10}}x}}{\left(1 + e^{\frac{3}{\sqrt{10}}x}\right)^{\frac{11}{3}}}\right)\frac{t^3}{3}$$

By setting $p = 1$ in Eq. (27) the solution of Eq. (24) can be obtained as $v = v_0 + v_1 + v_2 + v_3 + \cdots$ thus the solution of Eq. (24) can be written as:

$$v(x,t) = \left(\frac{1}{1 + e^{\frac{3}{\sqrt{10}}x}}\right)^{\frac{2}{3}} + \frac{7}{5}\left(\frac{e^{\frac{3}{\sqrt{10}}x}}{\left(1 + e^{\frac{3}{\sqrt{10}}x}\right)^{\frac{5}{3}}}\right) t + \frac{49}{50}\left(\frac{\left(2e^{\frac{3}{\sqrt{10}}x} - 3\right)e^{\frac{3}{\sqrt{10}}x}}{\left(1 + e^{\frac{3}{\sqrt{10}}x}\right)^{\frac{8}{3}}}\right)\frac{t^2}{2}$$

$$+ \frac{343}{1000}\left(\frac{\left(4\left(e^{\frac{3}{\sqrt{10}}x}\right)^2 - 27e^{\frac{3}{\sqrt{10}}x} + 9\right)e^{\frac{3}{\sqrt{10}}x}}{\left(1 + e^{\frac{3}{\sqrt{10}}x}\right)^{\frac{11}{3}}}\right)\frac{t^3}{3} + \cdots \qquad (30)$$

The Taylor series expansion for $\left( \frac{1}{2} tanh \left( -\frac{3}{2\sqrt{10}} \left( x - \frac{7}{\sqrt{10}} t \right) \right) + \frac{1}{2} \right)^{\frac{2}{3}}$ is written as:

$$\left( \frac{1}{2} tanh \left( -\frac{3}{2\sqrt{10}} \left( x - \frac{7}{\sqrt{10}} t \right) \right) + \frac{1}{2} \right)^{\frac{2}{3}}$$

$$= \left( \frac{1}{1 + e^{\frac{3}{\sqrt{10}}x}} \right)^{\frac{2}{3}} + \frac{7}{5} \left( \frac{e^{\frac{3}{\sqrt{10}}x}}{\left( 1 + e^{\frac{3}{\sqrt{10}}x} \right)^{\frac{5}{3}}} \right) t + \frac{49}{50} \left( \frac{\left( 2 e^{\frac{3}{\sqrt{10}}x} - 3 \right) e^{\frac{3}{\sqrt{10}}x}}{\left( 1 + e^{\frac{3}{\sqrt{10}}x} \right)^{\frac{8}{3}}} \right) \frac{t^2}{2}$$

$$+ \frac{343}{1000} \left( \frac{\left( 4 \left( e^{\frac{3}{\sqrt{10}}x} \right)^2 - 27 e^{\frac{3}{\sqrt{10}}x} + 9 \right) e^{\frac{3}{\sqrt{10}}x}}{\left( 1 + e^{\frac{3}{\sqrt{10}}x} \right)^{\frac{11}{3}}} \right) \frac{t^3}{3} + \cdots \qquad (31)$$

Comparing Eq. (31) with Eq. (30), thus Eq. (30) can be reduced to

$$v(x,t) = \left( \frac{1}{2} tanh \left( -\frac{3}{2\sqrt{10}} \left( x - \frac{7}{\sqrt{10}} t \right) \right) + \frac{1}{2} \right)^{\frac{2}{3}} \qquad (32)$$

This is the exact solution of the problem, Eq. (24). Table 3 shows the trend of rapid convergence of the results of $S_0(x,t) = v_0(x,t)$ to $S_5(x,t) = \sum_{i=0}^{5} v_i(x,t)$ using the HPM solution toward the exact solution. The maximum relative error of less than 0.007% is achieved in comparison to the exact solution as shown in table 3. We can conclude that the HPM is one most suitable and friendly method in solving the Newell-Whitehead-Segel equation.

Table 3 shows: the percentage of relative errors of the results of $S_0(x,t) = v_0(x,t)$ to $S_5(x,t) = \sum_{i=0}^{5} v_i(x,t)$ of the HPM solution of Eq. (24).

***Case IV:***

In this case we will examine the Newell-Whitehead-Segel equation for $a = 3, b = 4, k = 1, q = 3$,

$$\frac{\partial u}{\partial t} = \frac{\partial^2 u}{\partial x^2} + 3u - 4u^3, \qquad (33)$$

Subject to initial condition

$$u(x,0) = \sqrt{\frac{3}{4}} \frac{e^{\sqrt{6}x}}{e^{\sqrt{6}x} + e^{\frac{\sqrt{6}}{2}x}}, \qquad (34)$$

We construct a homotopy for Eq. (33) in the following form:

$$H(v,p) = (1-p) \left[ \frac{\partial v}{\partial t} - \frac{\partial u_0}{\partial t} \right] + p \left[ \frac{\partial v}{\partial t} - \frac{\partial^2 v}{\partial x^2} - 3v + 4v^3 \right] = 0. \qquad (35)$$

The solution of Eq. (33) can be written as a power series in $\boldsymbol{p}$ as:

$$v = v_0 + p\, v_1 + p^2 v_2 + \cdots \qquad (36)$$

Substituting Eq. (36) and Eq. (34) into Eq. (35) and equating the terms with identical powers of *p*:

$$p^0: \frac{\partial v_0}{\partial t} = \frac{\partial u_0}{\partial t}, \qquad v_0(x,0) = \sqrt{\frac{3}{4}} \frac{e^{\sqrt{6}x}}{e^{\sqrt{6}x} + e^{\frac{\sqrt{6}}{2}x}},$$

$$p^1: \frac{\partial v_1}{\partial t} = \frac{\partial^2 v_0}{\partial x^2} - \frac{\partial v_0}{\partial t} + v_0(3 - 4v_0{}^2) \qquad v_1(x,0) = 0,$$

$$p^2: \frac{\partial v_2}{\partial t} = \frac{\partial^2 v_1}{\partial x^2} - 12v_0{}^2 v_1 + 3v_1 \qquad v_2(x,0) = 0, \tag{37}$$

$$p^3: \frac{\partial v_3}{\partial t} = \frac{\partial^2 v_2}{\partial x^2} + 3v_2 - 12v_0 v_1^2 - 12 v_2 v_0^2 \qquad v_3(x,0) = 0.$$

**Table 3:**

| | | percentage of relative error (%RE) | | |
|---|---|---|---|---|
| | | $x = 1$ | $x = 1.5$ | $x = 1.8$ |
| $t = 0.1$ | $S_0(x,t)$ | 0.09323514238 | 0.1045372618 | 0.1099789256 |
| | $S_1(x,t)$ | 0.001725301473 | 0.003516790607 | 0.004498961488 |
| | $S_2(x,t)$ | 0.0002100221036 | 0.0001216772460 | 0.00005170881962 |
| | $S_3(x,t)$ | 0.00001173612982 | 0.00001539927439 | 0.00001472293353 |
| | $S_4(x,t)$ | 3.620778386 e-7 | 2.929797495 e-7 | 5.746720664 e-7 |
| | $S_5(x,t)$ | 5.414213656 e-8 | 3.614791446 e-8 | 1.434500497 e-8 |
| $t = .3$ | $S_0(x,t)$ | 0.2396054415 | 0.2702774140 | 0.2852064054 |
| | $S_1(x,t)$ | 0.009390569965 | 0.02330928345 | 0.03106724554 |
| | $S_2(x,t)$ | 0.005215736776 | 0.003375969333 | 0.001825357483 |
| | $S_3(x,t)$ | 0.0007262150389 | 0.001037577469 | 0.001023346798 |
| | $S_4(x,t)$ | 0.00009555513293 | 0.00004044412369 | 0.0001029640746 |
| | $S_5(x,t)$ | 0.00003280557979 | 0.00002473085989 | 0.00001198739229 |
| $t = .4$ | $S_0(x,t)$ | 0.2962707966 | 0.3351845099 | 0.3542406943 |
| | $S_1(x,t)$ | 0.01219208796 | 0.03518330953 | 0.04811460219 |
| | $S_2(x,t)$ | 0.01183961035 | 0.008037424274 | 0.004713588643 |
| | $S_3(x,t)$ | 0.001990820701 | 0.002987593293 | 0.002996129616 |
| | $S_4(x,t)$ | 0.0004128329726 | 0.0001164720814 | 0.0003682018367 |
| | $S_5(x,t)$ | 0.0001681121669 | 0.0001337455766 | 0.00006941923143 |

Using the Maple package to solve recursive sequences, Eq. (37), we obtain:

$$v_0(x,t) = \sqrt{\frac{3}{4}} \frac{e^{\sqrt{6}x}}{e^{\sqrt{6}x} + e^{\frac{\sqrt{6}}{2}x}},$$

$$v_1(x,t) = \frac{9}{2}\sqrt{\frac{3}{4}} \frac{e^{\sqrt{6}x} e^{\frac{\sqrt{6}}{2}x}}{\left(e^{\sqrt{6}x} + e^{\frac{\sqrt{6}}{2}x}\right)^2} t, \tag{38}$$

$$v_2(x,t) = \frac{81}{4}\sqrt{\frac{3}{4}} \frac{e^{\sqrt{6}x} e^{\frac{\sqrt{6}}{2}x} \left(-e^{\sqrt{6}x} + e^{\frac{\sqrt{6}}{2}x}\right)}{\left(e^{\sqrt{6}x} + e^{\frac{\sqrt{6}}{2}x}\right)^3} \frac{t^2}{2}$$

$$v_3(x,t) = \frac{243}{16}\sqrt{\frac{3}{4}}\frac{e^{\sqrt{6}x}e^{\frac{\sqrt{6}}{2}x}\left(-4e^{\sqrt{6}x}e^{\frac{\sqrt{6}}{2}x}+\left(e^{\sqrt{6}x}\right)^2+\left(e^{\frac{\sqrt{6}}{2}x}\right)^2\right)}{\left(e^{\sqrt{6}x}+e^{\frac{\sqrt{6}}{2}x}\right)^4}\frac{t^3}{3}$$

By setting $p = 1$ in Eq. (36) the solution of Eq. (33) can be obtained as $v = v_0 + v_1 + v_2 + v_3 + \dots$

Thus the solution of Eq. (33) can be written as,

$$v(x,t) =$$

$$\sqrt{\frac{3}{4}}\frac{e^{\sqrt{6}x}}{e^{\sqrt{6}x}+e^{\frac{\sqrt{6}}{2}x}} + \frac{9}{2}\sqrt{\frac{3}{4}}\frac{e^{\sqrt{6}x}e^{\frac{\sqrt{6}}{2}x}}{\left(e^{\sqrt{6}x}+e^{\frac{\sqrt{6}}{2}x}\right)^2}t + \frac{81}{4}\sqrt{\frac{3}{4}}\frac{e^{\sqrt{6}x}e^{\frac{\sqrt{6}}{2}x}\left(-e^{\sqrt{6}x}+e^{\frac{\sqrt{6}}{2}x}\right)}{\left(e^{\sqrt{6}x}+e^{\frac{\sqrt{6}}{2}x}\right)^3}\frac{t^2}{2} +$$

$$\frac{243}{16}\sqrt{\frac{3}{4}}\frac{e^{\sqrt{6}x}e^{\frac{\sqrt{6}}{2}x}\left(-4e^{\sqrt{6}x}e^{\frac{\sqrt{6}}{2}x}+\left(e^{\sqrt{6}x}\right)^2+\left(e^{\frac{\sqrt{6}}{2}x}\right)^2\right)}{\left(e^{\sqrt{6}x}+e^{\frac{\sqrt{6}}{2}x}\right)^4}\frac{t^3}{3} + \dots \qquad (39)$$

The Taylor series expansion for $\sqrt{\frac{3}{4}}\frac{e^{\sqrt{6}x}}{e^{\sqrt{6}x}+e^{\left(\frac{\sqrt{6}}{2}x-\frac{9}{2}t\right)}}$ is written as

$$\sqrt{\frac{3}{4}}\frac{e^{\sqrt{6}x}}{e^{\sqrt{6}x}+e^{\left(\frac{\sqrt{6}}{2}x-\frac{9}{2}t\right)}}$$

$$= \sqrt{\frac{3}{4}}\frac{e^{\sqrt{6}x}}{e^{\sqrt{6}x}+e^{\frac{\sqrt{6}}{2}x}} + \frac{9}{2}\sqrt{\frac{3}{4}}\frac{e^{\sqrt{6}x}e^{\frac{\sqrt{6}}{2}x}}{\left(e^{\sqrt{6}x}+e^{\frac{\sqrt{6}}{2}x}\right)^2}t$$

$$+ \frac{81}{4}\sqrt{\frac{3}{4}}\frac{e^{\sqrt{6}x}e^{\frac{\sqrt{6}}{2}x}\left(-e^{\sqrt{6}x}+e^{\frac{\sqrt{6}}{2}x}\right)}{\left(e^{\sqrt{6}x}+e^{\frac{\sqrt{6}}{2}x}\right)^3}\frac{t^2}{2}$$

$$+ \frac{243}{16}\sqrt{\frac{3}{4}}\frac{e^{\sqrt{6}x}e^{\frac{\sqrt{6}}{2}x}\left(-4e^{\sqrt{6}x}e^{\frac{\sqrt{6}}{2}x}+\left(e^{\sqrt{6}x}\right)^2+\left(e^{\frac{\sqrt{6}}{2}x}\right)^2\right)}{\left(e^{\sqrt{6}x}+e^{\frac{\sqrt{6}}{2}x}\right)^4}\frac{t^3}{3} + \dots \qquad (40)$$

By substituting Eq. (40) into Eq. (39), the Eq. (39) can be reduced to:

$$v(x,t) = \sqrt{\frac{3}{4}}\frac{e^{\sqrt{6}x}}{e^{\sqrt{6}x}+e^{\left(\frac{\sqrt{6}}{2}x-\frac{9}{2}t\right)}} \qquad (41)$$

This is the exact solution of the problem, Eq. (33). Table 4 shows the trend of rapid convergence of the results of $S_0(x,t) = v_0(x,t)$ to $S_5(x,t) = \sum_{i=0}^{5} v_i(x,t)$ using the HPM solution toward the exact solution. The maximum relative error of less than 0.037% is achieved in comparison to the exact solution as shown in table 4.

**Table 4:**

| | | percentage of relative error (%RE) | | |
|---|---|---|---|---|
| | | $x = .2$ | $x = .4$ | $x = .8$ |
| $t = 0.1$ | $S_0(x,t)$ | 0.1591055626 | 0.1376714465 | 0.09890297594 |
| | $S_1(x,t)$ | 0.007038530707 | 0.009754764500 | 0.01176933513 |
| | $S_2(x,t)$ | 0.002482893305 | 0.001788283111 | 0.0004607861238 |
| | $S_3(x,t)$ | 0.0001958736131 | 0.0002690488417 | 0.0002512962487 |
| | $S_4(x,t)$ | 0.00004554598541 | 0.00002344226507 | 0.00001229642845 |
| | $S_5(x,t)$ | 0.000004998678372 | 0.000006399733729 | 0.000003681928211 |
| $t = .15$ | $S_0(x,t)$ | 0.2155132753 | 0.1864801193 | 0.1339670591 |
| | $S_1(x,t)$ | 0.01698529379 | 0.02214248709 | 0.02558157620 |
| | $S_2(x,t)$ | 0.007422698254 | 0.005232453593 | 0.001127443849 |
| | $S_3(x,t)$ | 0.001011675135 | 0.001318032849 | 0.001182316216 |
| | $S_4(x,t)$ | 0.0003016921205 | 0.0001450263132 | 0.0001001951821 |
| | $S_5(x,t)$ | 0.00005638430243 | 0.00006875985734 | 0.00003732434558 |
| $t = .2$ | $S_0(x,t)$ | 0.2605557486 | 0.2254546357 | 0.1619662983 |
| | $S_1(x,t)$ | 0.03164334146 | 0.03938244479 | 0.04388750886 |
| | $S_2(x,t)$ | 0.01561926168 | 0.01076040185 | 0.001819030418 |
| | $S_3(x,t)$ | 0.003225422844 | 0.004022798070 | 0.003478947590 |
| | $S_4(x,t)$ | 0.001110362592 | 0.0004931283715 | 0.0004433744920 |
| | $S_5(x,t)$ | 0.0003119309814 | 0.0003646043520 | 0.0001870028440 |

**Table 4** shows: the percentage of relative errors of the results of $S_0(x,t) = v_0(x,t)$ to $S_5(x,t) = \sum_{i=0}^{5} v_i(x,t)$ of the HPM solution of Eq. (33).

### 4. Conclusion:

In the present work the exact solution of the Newell-Whitehead-Segel nonlinear diffusion equation is obtained using the HPM. The validity and effectiveness of the HPM is shown by solving four non-homogenous non-linear differential equations and the very rapid convergence to the exact solutions is demonstrated numerically. The trend of rapid and monotonic convergence of the solution toward the exact solution is clearly shown by obtaining the relative error in compared to the exact solution. The rapid convergence towards the exact solutions of HPM indicates that, using the HPM to solve the non-linear differential equations, a reasonable less amount of computational work with acceptable accuracy may be sufficient. Moreover it can be concluded that the HPM is a very powerful and efficient technique which can construct the exact solution of nonlinear differential equations.